\documentclass{article}
\usepackage{amsthm,amsmath,amssymb}
\theoremstyle{plain}
\newtheorem{theorem}{Theorem}[section]

\newtheorem{proposition}[theorem]{Proposition}
\newtheorem{corollary}[theorem]{Corollary}

\theoremstyle{definition}
\newtheorem{definition}[theorem]{Definition}

\theoremstyle{remark}

\title{Improving an algorithm to solve Multiple Simultaneous
Conjugacy Problems in braid groups}
\author{Juan Gonz\'alez-Meneses}
\date{November 2002.}

\begin{document}
\maketitle

\begin{abstract}
 There are recent cryptographic protocols that are based on Multiple Simultaneous Conjugacy Problems
in braid groups. We improve an algorithm, due to Sang Jin Lee and Eonkyung Lee, to solve
these problems, by applying a method developed by the author and Nuno Franco,
originally intended to solve the Conjugacy Search Problem in braid groups.
\end{abstract}

\section{Introduction}
 In~\cite{leelee}, Sang Jin Lee and Eonkyung Lee give an algorithm to solve the following problem,
that they call Multiple Simultaneous Conjugacy Problem (MSCP), in the braid group $B_n$: given
the $r$-tuples $(a_1,\ldots, a_r)$ and  $(x^{-1}a_1x,\ldots, x^{-1}a_rx)$ in $B_n$, find the
conjugator $x$.

 This problem has been proposed for cryptographical applications: There is a Key
Agreement Protocol proposed by Anshell, Anshell and Goldfeld in \cite{AAG}, improved by the same authors
and Fisher in~\cite{AAFG}, which is based on the difficulty to solve a MSCP in some groups.
Braid groups have been proposed as a good choice. There have been different attacks to this
cryptosystem, namely length-based attacks (\cite{HT}, \cite{bar-ilan}), linear algebraic ones
(\cite{leelee}, \cite{H}) and others (\cite{rainer}). But the algorithm we describe in this paper can
be thought of as a direct attack to the base problem of the protocol.

 We will assume that the reader is familiar with the basic notions in braid theory, which can be found
in \cite{B} or \cite{MK}. It is also desirable to know the work in \cite{garside}, \cite{EM} and \cite{T}.

 Recall that, given a braid $a\in B_n$, the integer $\inf(a)$ is the biggest $k\in \mathbb{Z}$ such that
$a=\Delta^{k}p$, where $\Delta$ is the usual Garside element (half twist of all the strands) and $p$ is
a positive braid (all its crossings are positive).

 The algorithm in~\cite{leelee} works as follows: First they define, for every $r$-tuple of braids,
$\alpha=(a_1,\ldots, a_r) \in (B_n)^r$, the set $C^{\inf}(\alpha)$ consisting of all
$\beta=(b_1,\ldots,b_r)\in (B_n)^r$ such that $\inf(b_i)\geq \inf(a_i)$ for all $i$ and there exists some
$\omega\in B_n$ satisfying $b_i=\omega^{-1}a_i\omega$ for all $i$ simultaneously (that is,
$\beta=\omega^{-1} \alpha \omega$).
Then they prove the following result:

\begin{theorem}{\em \cite{leelee}}\label{T:beta'}
Let $\alpha=(a_1,\ldots, a_r)$ and $\beta=(b_1,\ldots, b_r)$ be an instance of a MSCP in $B_n$,
and $x$ a positive solution. Then one can compute a positive braid $x_0$ and a $r$-tuple
$\beta'=(b_1',\ldots, b_r')\in C^{\inf}(\alpha)$ such that $b_i'=x_0b_ix_0^{-1}$
for all $i$, in time proportional to
$$
 n(\log n)|x|\left( |x|+\sum_{i=1}^r{(|a_i|+|b_i|)}\right),
$$
where $|\cdot|$ denotes word length in generators. Moreover $x=x_1x_0$ for some positive braid
$x_1$.
\end{theorem}

Here $C^{\inf}(\alpha)$ plays the role of the {\em Summit Set} defined
in~\cite{garside} to solve the conjugacy problem in $B_n$, in the sense that it satisfies
the following result:

\begin{theorem}\label{T:chain}{\em \cite{leelee}}
Given $\beta\in C^{\inf}(\alpha)$, there exists a chain of elements
$\alpha=\alpha_1,\alpha_2,\ldots,\alpha_{k+1}=\beta$
in $C^{\inf}(\alpha)$, where successive elements are simultaneously conjugated by
a permutation braid. In other words, there exist permutation braids $s_1,\ldots,s_k$ such that
$s_j^{-1} \alpha_j s_j = \alpha_{j+1}$ for every $j=1,\ldots,k$.
\end{theorem}

Therefore, by classical methods (see~\cite{garside}), one can use these two results to solve any
MSCP in finite time. Nevertheless, this classical approach gives a computational complexity which is
exponential with respect to the braid index $n$, and involves the cardinality $N$ of the set
$C^{\inf}(\alpha)$.

S. J. Lee and E. Lee expect in~\cite{leelee} that one can apply the methods in~\cite{francogm} to
this algorithm, so that the computational complexity becomes a polynomial in $(n,r,l,N)$, where
$l$ is the maximal word-length of the $a_i$'s and $b_i$'s. Here we show that this is the case. More precisely,
we show:

\begin{theorem}\label{T:complexity}
 Let $\alpha=(a_1,\ldots,a_r)\in(B_n)^r$ and let $\beta=(b_1,\ldots,b_r)\in C^{\inf}(\alpha)$.
Let $l$ be the maximal word length of the $a_i$'s and $b_i$'s, and let $N$ be the number of elements in
$C^{\inf}(\alpha)$. Then one can compute a braid $x\in B_n$ such that $x^{-1}\alpha x=\beta$ in time
$O(Nrl^2n^3)$.
\end{theorem}

\section{Minimal simple elements for MSCP}

 Let us consider the {\em Artin monoid of positive braids}, $B_n^+$. We can define a {\em prefix order}
on its elements, $\prec$, as follows: for $a,b\in B_n^+$, $a\prec b$ if and only if there exists $c\in B_n^+$
such that $ac=b$. We will say that $a$ is a {\em prefix} (or a {\em divisor}) of $b$, or that $b$ is
divisible by $a$. This is a partial order on $B_n^+$, with some nice properties: For every $u, v\in B_n^+$
there exists their {\em least common multiple}, $u\vee v$, and their {\em greatest common divisor},
$u\wedge v$. There also exists an element $\Delta$ (which is represented by a half twist of all the
strands) which, together with the above partial order, endows $B_n^+$ with a structure of {\em Garside
monoid}, so $B_n$ is a {\em Garside group} (cf. \cite{DP} \cite{D}).

 The {\em permutation braids}, also called {\em simple elements}, are the prefixes (or divisors)
of $\Delta$. We denote by $S$ the set of simple elements. In $B_n^+$ there are $n!$ simple elements.

 The algorithm used in~\cite{leelee} to solve a MSCP goes as follows: given
$\alpha, \beta \in (B_n)^r$ conjugated, one computes $\beta'\in C^{\inf}(\alpha)$ as in Theorem~\ref{T:beta'}.
Then one must construct the whole $C^{\inf}(\alpha)$ using the method by Garside:
Conjugate $\alpha$ by all simple elements. If new elements in $C^{\inf}(\alpha)$ are obtained, conjugate
each one of them by all simple elements. Continue until no new elements appear. At that point, by
Theorem~\ref{T:chain}, we will have computed the whole $C^{\inf}(\alpha)$ and moreover, we will know a
chain going from $\alpha$ to any other element in $C^{\inf}(\alpha)$, as in Theorem~\ref{T:chain}.
Hence, the chain associated to $\beta'$, together with the element $x_0$ in Theorem~\ref{T:beta'} will
give us the solution to the MSCP.

 One of the main problems of this algorithm is the size of $S$.
For every element in $C^{\inf}(\alpha)$ one must compute $n!$ conjugations!
The idea in~\cite{francogm} is to consider very small subsets of $S$, which can be fastly computed,
satisfying some suitable properties that allow the classical algorithm to work with them, instead of the
whole $S$. The general method to compute these small subsets is the following.

 Let ${\cal P}$ be a property for simple elements, and let $S_{\cal P}$ be the set of simple elements
satisfying ${\cal P}$. Then $\min(S_{\cal P})$ is defined as the set of minimal elements (with
respect to $\prec$) in $S_{\cal P}$. We must then define some suitable properties.

  Let $J=(j_1,\ldots,j_r)\in \mathbb{Z}^r$ and let $C_J$ be the set of $r$-tuples
$\delta=(d_1,\ldots, d_r)\in (B_n)^r$ such that $\inf(d_i)\geq j_i$ for all $i$.

\begin{definition}
 Let $J=(j_1,\ldots,j_r)\in \mathbb{Z}^r$ and let $\delta=(d_1,\ldots, d_r)\in C_J$. We say that a
simple element $s$ satisfies the property ${\cal P}(\delta,J)$ if $s^{-1}\delta s \in C_J$. In other words,
if $\inf(s^{-1}d_is)\geq j_i$ for all $i$.
\end{definition}

 Now consider the subsets $S_{\delta,J}=\min(S_{{\cal P}(\delta,J)})\subset S$, where
$\delta\in C_J$. These are the small subsets of $S$ we were talking about.
We can use them to solve a MSCP by means of the following result:

\begin{proposition}\label{P:chain}
 Given $\alpha=(a_1,\ldots,a_r)\in (B_n)^r$, let $J=(\inf(a_1),\ldots,\inf(a_r))\in \mathbb{Z}^r$.
For every $\beta\in C^{\inf}(\alpha)$, there exists a chain
$\alpha=\alpha_1,\alpha_2,\ldots,\alpha_{k+1}=\beta$ in $C^{\inf}(\alpha)$,
where for $j=1,\ldots,k$, $\;\alpha_j$ is conjugated to $\alpha_{j+1}$ by a simple element
$s_j\in S_{\alpha_j,J}$. That is, $s_j^{-1} \alpha_j s_j = \alpha_{j+1}$ and $s_j$ is minimal among
the simple elements conjugating $\alpha_j$ to an element in $C_J$.
\end{proposition}

\begin{proof}
  This result is analogous to Proposition~4.10 in \cite{francogm}. It suffices to take the chain given
in Theorem~\ref{T:chain} and decompose every simple element into minimal ones. We notice that we
obtain a chain of elements in $C_J$, but since all these elements are conjugated to $\alpha$, they all
belong to $C^{\inf}(\alpha)$.
\end{proof}

\section{Size of $\mathbf{S_{\delta,J}}$}

 In this section we will show that the cardinal of $S_{\delta,J}$, for every $J\in\mathbb{Z}^r$ and
every $\delta\in C_J$, is always smaller that $n$. Hence, if we know how to compute it fastly, we
will improve considerably the speed of the algorithm by Lee and Lee (recall that $\#(S)=n!$).
We will need the following results:

\begin{proposition}\label{P:closed_gcd->small}{\em \cite{francogm}}
If a property ${\cal P}$ is closed under gcd (i.e., if $s_1,s_2\in S_{\cal P}$ implies
$s_1\wedge s_2\in S_{\cal P}$) then $\#(\min(S_{\cal P}))\leq n-1$. \hfill \fbox{}
\end{proposition}

\begin{proposition}\label{closed_gcd} For every $J\in \mathbb{Z}^r$ and every $\delta\in C_J$,
the property ${\cal P}(\delta,J)$ is closed under gcd.
\end{proposition}

\begin{proof}
 Suppose that $s_1,s_2\in S_{{\cal P}(\delta,J)}$, that is, for every $i=1,\ldots,r$,
$\inf(s_1^{-1}d_is_1)\geq j_i$ and $\inf(s_2^{-1}d_is_2)\geq j_i$. Since $\delta\in C_J$
one has $d_i=\Delta^{j_i}p_i$ for some positive braid $p_i$. Then
$$
   s_1^{-1}d_is_1 =  s_1^{-1}\Delta^{j_i}p_is_1 = \Delta^{j_i}\tau^{j_i}(s_1^{-1})p_is_1,
$$
where $\tau$ is the inner automorphism of $B_n$ which consists on conjugation by $\Delta$.
Hence, $\inf(s_1^{-1}d_is_1)\geq j_i$ means that $\tau^{j_i}(s_1^{-1})p_is_1$ is positive,
or in other words: $\tau^{j_i}(s_1)\prec p_is_1$. In the same way one has
$\tau^{j_i}(s_2)\prec p_is_2$ for all $i$. We must therefore show that, for $i=1,\ldots,r$,
one has $\tau^{j_i}(s)\prec p_is$, where $s=s_1\wedge s_2$.

 Since $\tau$ is a homomorphism that preserves the prefix order, then
$\tau^{j_i}(s_1)\wedge \tau^{j_i}(s_2) = \tau^{j_i}(s_1\wedge s_2)=\tau^{j_i}(s)$.
This implies $\tau^{j_i}(s)\prec p_i s_1$ and $\tau^{j_i}(s)\prec p_i s_2$, hence
$\tau^{j_i}(s)\prec (p_is_1)\wedge (p_i s_2)= p_i(s_1\wedge s_2)= p_i s$, as we wanted to show.
\end{proof}

\begin{corollary}\label{C:small}
 For every $J\in \mathbb{Z}^r$ and every $\delta\in C_J$, the set
$S_{\beta,J}=\min(S_{{\cal P}(\beta,J)})$ has at most $n-1$ elements. \hfill \fbox{}
\end{corollary}

\section{How to compute $\mathbf{S_{\delta,J}}$}

 We will finally present an algorithm that computes $S_{\delta,J}$, given $J\in \mathbb{Z}^r$ and
$\delta\in C_J$. This algorithm will have complexity $O(rl^2n^3)$. Hence, in the algorithm by Lee and Lee,
we no longer need to conjugate every $\delta\in C^{\inf}(\alpha)$ by all simple elements
(n! conjugations); we can compute $S_{\delta,J}$ and then we do no more than $n-1$ conjugations.

 We first need to be more precise about the work in~\cite{francogm}. We saw in
Proposition~\ref{P:closed_gcd->small} that $\min(S_{\cal P})$ has at most $n-1$ elements; but be can
actually say more: for every generator $\sigma_i$, there is exactly one element
$r_i\in \min(S_{\cal P})$ such that $\sigma_i\prec r_i$. It can happen, however, that $r_i=r_j$ for some
$i\neq j$. Anyway, in order to compute $\min(S_{\cal P})$ (in our particular case $S_{\delta,J}$), we
just need to compute $r_i$ for $i=1,\ldots,n-1$.

 It is also given in~\cite{francogm} a method to compute the least common multiple $s\vee p$ of a simple
element $s$ and a positive braid $p$. More precisely, the algorithm given in~\cite{francogm} computes
a simple element $s'$ such that $ps'=s\vee p$. This takes time $O(l^2n\log n)$, where $l$ is the word
length of $p$, and $n$ is the number of strands. Notice that, in terms of theoretical complexity, this
algorithm is equivalent to the computation a normal form (cf. \cite{T}). Furthermore, it is also shown
in~\cite{francogm} that if $p$ is given in left normal form, then the complexity becomes $O(ln\log n)$.

 So let us suppose that we are given $J=(j_1,\ldots, j_r)\in \mathbb{Z}^r$ and $\delta=(d_1,\ldots,d_r)\in C_J$, and we
want to compute $S_{\delta,J}$. As we said before, we just need to compute $r_i$ for every
$i=1,\ldots,n-1$, where, in this case, $r_i$ is the minimal simple element which is divisible by
$\sigma_i$ and conjugates $\delta$ to an element in $C_J$. We propose the following algorithm:

\vspace{.3cm}
\noindent \underline{\bf Algorithm to compute $\mathbf{r_i}$.}
\begin{enumerate}

 \item Let $D\subset \{1,\ldots,r\}$ consisting of those $t$ such that $\inf(d_t)=j_t$.

 \item For every $t\in D$, compute $p_t$ such that $d_t=\Delta^{j_t}p_t$.

 \item Let $s=\sigma_i$.

 \item If $\tau^{j_t}(s)\prec p_ts$ for every $t\in D$, then return $s$. Stop.

 \item Take $m\in D$ such that $\tau^{j_m}(s)\not\prec p_ms$.

 \item Compute $s'$ such that $(p_m s)s'=\tau^{j_m}(s)\vee p_m s$.

 \item Let $s=s s'$ and go to step 4.

\end{enumerate}

\begin{proposition}
 Given $J=(j_1,\ldots, j_r)\in \mathbb{Z}^r$, $\delta=(d_1,\ldots,d_r)\in C_J$ and
$i\in \{1,\ldots,n-1\}$, the above algorithm computes $r_i$, the minimal simple element
which is divisible by $\sigma_i$ and conjugates $\delta$ to an element in $C_J$.
\end{proposition}

\begin{proof}
 The algorithm starts by considering just those $d_t$ whose infimum is exactly $j_t$. This is due
to the following fact: If we can write $d_t=\Delta^{k}p_t$ where $k> j_t$ and $p_t$ is a positive braid,
then for every simple element $s$ we will have:
$$
s^{-1}d_ts=s^{-1}\Delta^k p_t s = \Delta^k \tau^k(s^{-1}) p_t s =
\Delta^{k-1} (\Delta \tau^k(s^{-1}))p_t s.
$$
But $\tau^k(s)$ is a simple element, so $\Delta \tau^k(s^{-1})$ is a positive braid, hence
the infimum of $s^{-1}d_ts$ is at least $k-1\geq j_t$. Therefore, we just need to care about those
$d_t$ where $t\in D$.

 For every $t\in D$ one has $d_t=\Delta^{j_t}p_t$, where $p_t$ is a positive braid. These elements
$p_t$ are computed in Step 2 just by computing the left normal form of $d_t$.

 We want to find $r_i$, and we know that $\sigma_i\prec r_i$. In the algorithm, the simple element $s$
will be the possible value of $r_i$. At every iteration of the loop in steps 4-7, we start with a
simple element $s$ such that $\sigma_i \prec  s \prec r_i$, and we check if $s=r_i$. If it is not,
we multiply $s$ by some suitable simple element $s'$, and we start again. We must show that this
makes sense.

 At Step 3 we set $s=\sigma_i$, so we are sure that $\sigma_i\prec s \prec r_i$. Then we start the loop.
In order to decide if $s=r_i$, we must check if $\inf(s^{-1}d_ts)\geq j_t$ for all $t\in D$.
But, in the same way as above, one has $s^{-1}d_ts= \Delta^{j_t}\tau^{j_t}(s^{-1}) p_t s$, so
$\inf(s^{-1}d_ts)\geq j_t$ if and only if $\tau^{j_t}(s^{-1}) p_t s$ is a positive braid, or in other
words, if $\tau^{j_t}(s)\prec p_ts$. This is what is checked at Step 4.

 If Step 4 determined that $s\neq r_t$, we must have found some $m\in D$ such that
$\tau^{j_m}(s)\not\prec p_ms$. Step 5 just takes one of these values.

Now it comes the main step: We know that $s\prec r_i$, so $r_i=s\widehat{s}$ for some simple element
$\widehat{s}$. Moreover, $\inf(r_i^{-1}d_mr_i)\geq j_t$ so one has $\tau^{j_m}(r_i)\prec p_mr_i$.
Hence,  $\tau^{j_m}(s) \prec \tau^{j_m}(s)\tau^{j_m}(\widehat{s})= \tau^{j_m}(r_i) \prec p_mr_i$
while on the other hand $p_m s\prec p_m s \widehat{s} =p_m r_i$. Therefore, the least common
multiple $\tau^{j_m}(s)\vee p_m s$ must also divide $p_m r_i$. Step~6 computes this lcm. Actually,
it computes $s'$ such that $\tau^{j_m}(s)\vee p_m s=(p_ms)s'$. But since this divides $p_mr_i$,
we finally obtain that $ss'\prec r_i$.

 We must remark two facts: First, $ss'$ is always a simple element, since it divides the simple element
$r_i$. Second, $s'$ cannot be trivial, since otherwise we would have $\tau^{j_m}(s)\vee p_m s= p_m s$,
implying $\tau^{j_m}(s)\prec p_m s$, which gives a contradiction with the choice of $m$.
Therefore, $ss'$ is strictly greater than $s$, but still a divisor of $r_i$, so in Step~7 we set
$s=ss'$, and start the loop again. This cannot run forever since the word length of $s$ is increased at
every iteration, so the maximal number of iterations is $\frac{n(n-1)}{2}$ (the word length of $\Delta$).

Therefore, at a certain iteration, we will obtain $s=r_i$, and the algorithm stops at Step~4 giving the
correct output.
\end{proof}

\section{Theoretical complexity}

 The algorithm we presented in this paper is exactly as the one in~\cite{leelee} except for the
computations of $S_{\delta,J}$, for every $\delta\in C^{\inf}(\alpha)$. The main step is the
computation of $r_i$ given by the algorithm in the previous section. So we start by studying
the complexity of this computation:

\begin{proposition}
  Given $J=(j_1,\ldots, j_r)\in \mathbb{Z}^r$, $\delta=(d_1,\ldots,d_r)\in C_J$ and
$i\in \{1,\ldots,n-1\}$, one can compute $r_i$ (the minimal simple element
which is divisible by $\sigma_i$ and conjugates $\delta$ to an element in $C_J$) in time
$O(rl^2n^2)$ where $l$ is the maximal word-length of the $d_i$'s.
\end{proposition}

\begin{proof}
 We need to study the complexity of the algorithm in the previous section.
First, Step~1 can be performed by computing the left normal form of every $d_t$. Every normal
form takes time $O(l^2n\log n)$, so Step one can be done in time $O(rl^2n\log n)$.

 The requirements of Step~2 can be achieved while doing Step~1: if some $d_t$ has infimum
 $j_t$, we keep the value of $p_t$. Hence Step~2 is negligible, as well as Step~3.

 Now we start a loop in Steps~4-7, which has at most $\frac{n(n-1)}{2}$ iterations, as we saw above.
The only non-negligible steps are Steps~4 and 6. In Step~4, for every $t\in D$ we must compute
$\tau^{j_t}(s)$, which can be done in linear time on the word size of $s$ (at most $\frac{n(n-1)}{2}$),
and then we must compute the left normal form of $p_t s$ taking time $O(ln\log n)$ (notice that $p_t$
is already in left normal form). After performing these computations, to check if
$\tau^{j_t}(s)\prec p_t s$ is $O(n\log n)$ (cf~\cite{T}). Hence Step~4 takes time $O(rln^2)$.
On the other hand, Step~6 can be done in time $O(ln\log n)$ by~\cite{francogm}.
Therefore, each iteration of the loop takes time $O(rln^2)$.

 Now we could say that, since there are at most $\frac{n(n-1)}{2}$ iterations, all of them can be computed
in time $O(rln^4)$. But we can do better than that: The different values of $s$ in the successive
iterations form an ascending chain of simple elements. Hence, the total number of computations performed in
all the iterations is the same as if it were just one iteration, with the maximum value of $s$
(see~\cite{T}). Therefore, the whole loop can be done in time $O(rln^2)$, and the whole algorithm takes
time $O(rl^2n^2)$.
\end{proof}

 We can now apply this result to measure our contribution to the algorithm in~\cite{leelee}:

\vspace{.3cm}
\noindent
{\it Proof of Theorem~\ref{T:complexity}.}
 One just need to apply the classical algorithm by Garside, together with the results given in
Proposition~\ref{P:chain} and Corollary~\ref{C:small}. To be more precise, let
$J=(\inf(a_1),\ldots,\inf(a_r))\in \mathbb{Z}^r$. For every element
$\delta\in C^{\inf}(\alpha)$ (there are $N$ elements) one must compute $S_{\delta,J}$. This takes
time $O(rl^2n^2)$ for every element, by the above result. Since there are at most $n-1$ elements,
it takes time $O(rl^2n^3)$. Then one must conjugate $\delta$ by all the elements
in $S_{\delta,J}$ (at most $n-1$), so we do at most $n-1$ conjugations by
simple elements, each one taking time $O(ln\log n)$ since $\delta$ is already in left normal form.

The algorithm stops when we find $\beta$. So, in the worst case, the complexity of the whole computation
is $O(Nrl^2n^3)$, as we wanted to show.  \hfill $\framebox{}$

\section{Final remarks}

In this paper we have improved the algorithm in~\cite{leelee} to solve a MSCP. More precisely, we have
improved a particular case of a MSCP, when the conjugate elements $\alpha$ and $\beta$ are such that
$\beta\in C^{\inf}(\alpha)$.

It is shown in~\cite{leelee} how to transform the general situation into this particular case
(see Theorem~\ref{T:beta'}), but the complexity of this step depends on the size of the solution!
Therefore, using this method we do not have an upper bound for the complexity of the general case,
in terms of the input data. Nevertheless, if our interest is to attack the cryptosystem in~\cite{AAG},
where the secret key is the solution to the MSCP, then the complexity given in Theorem~\ref{T:beta'},
to transform the general case into this particular case, yields a very efficient running time.

 Nevertheless, if one dislikes to measure the complexity in terms of the length of the solution, one
can do the following: given two conjugate elements $\alpha=(a_1,\ldots,a_r)$ and $\beta=(b_1,\ldots,b_r)$
in $(B_n)^r$, let $J=(j_1,\ldots,j_r)\in \mathbb{Z}^r$ where $j_i=\min(\inf(a_i),\inf(b_i))$. Then
one has $\alpha,\beta\in C_J$. Now define
$C^{\inf}(\alpha,\beta)$ as the set of $\delta\in C_J$ conjugate to $\alpha$ (thus to $\beta$).
Then all the above results can be applied to $C^{\inf}(\alpha,\beta)$, so we do not need to pass through
Theorem~\ref{T:beta'}. That is, we have:

\begin{theorem}\label{T:final}
    Let $\alpha=(a_1,\ldots,a_r)$ and $\beta=(b_1,\ldots,b_r)$ in $(B_n)^r$.
Let $l$ be the maximal word length of the $a_i$'s and $b_i$'s, and let $M$ be the number of elements in
$C^{\inf}(\alpha,\beta)$. Then one can compute a braid $x\in B_n$ such that $x^{-1}\alpha x=\beta$ in time
$O(Mrl^2n^3)$.
\end{theorem}

 Anyway, we do not think that this is the better way to proceed, since  $C^{\inf}(\alpha,\beta)$
will be, in general, much bigger than $C^{\inf}(\alpha)$, so one should try first to raise the infimum
of the entries of $\alpha$ and $\beta$, before starting to construct the whole $C^{\inf}(\alpha,\beta)$.

 On the other hand, the complexity given in Theorems~\ref{T:complexity} and \ref{T:final} may lead to
confusion, since one may think that we solved the MSCP in polynomial time. This is not true, since the
factors $N$ and $M$ (the size of $C^{\inf}(\alpha)$ and $C^{\inf}(\alpha,\beta)$) may not be a polynomial
in $(n,r,l)$ (there is no known bounds for $N$ or $M$ in terms of $(n,l,r)$). All we can say by now is
that $N$ and $M$ get smaller as $r$ grows, so it seems that MSCP's are simpler than usual conjugacy problems in
braid groups (see the discussion in~\cite{leelee} about the size of $N$).

 Finally, the algorithm in this paper works not only for braid groups, but for a larger
class of groups, called {\em Garside} groups (see~\cite{DP}, \cite{D} and \cite{francogm}), that
share with braid groups the existence of simple elements and their basic properties.
It can also be applied to other Garside structures in braid groups, as the one obtained from the
presentation by Birman, Ko and Lee in~\cite{BKL}.

\vspace{.3cm}
\noindent {\footnotesize
{\bf Juan Gonz\'alez-Meneses:} \\
Dep. Matem\'atica Aplicada I, ETS Arquitectura, Univ. de Sevilla,
Av. Reina Mercedes 2, 41012-Sevilla (SPAIN).
\\ E-mail: {\em meneses@us.es}}


\begin{thebibliography}{99}



\bibitem{AAFG} I. Anshel, M. Anshel, B. Fisher and D. Goldfeld, {\it New Key Agreement Protocols in Braid Group
Cryptography}. Topics in Cryptology--CT-RSA 2001 (San Francisco, CA), 13-27, Lecture Notes in Comput. Sci., {\bf 2020},
Springer, Berlin, 2001.


\bibitem{AAG} I. Anshel, M. Anshel and D. Goldfeld, {\it An algebraic method
for public-key cryptography}.  Math. Res. Lett. {\bf 6}, No. 3-4 (1999), 287-291.

\bibitem{B}
 J. Birman, ``Braids, links and mapping class groups".
Princeton University Press, Princeton, 1974.

\bibitem{BKL} J. Birman, K. H. Ko and S. J. Lee, {\it A new approach to
the word and conjugacy problems in the braid groups}, Adv. Math. {\bf 139},
No. 2 (1998), 322-353.


\bibitem{D} P. Dehornoy, {\it Groupes de Garside},  Ann. Scient. \'Ec. Norm. Sup., $4^e$~s\'erie, t. 35,
2002, 267-306.


\bibitem{DP} P. Dehornoy and L. Paris, {\it Gaussian groups and Garside
groups, two generalizations of Artin groups}, Proc. London Math. Soc. {\bf 79},
No. 3 (1999), 569-604.


\bibitem{EM} E. A. Elrifai, H. R. Morton, Algorithms for positive
braids, {\it Quart. J. Math. Oxford} {\bf 45} (1994), 479-497.

\bibitem{francogm} N. Franco and J. Gonz\'alez-Meneses, {\it Conjugacy problem for braid groups and
Garside groups}. To appear in Journal of Algebra. Available at {\tt www.arxiv.org/math.GT/0112310}

\bibitem{bar-ilan} D. Garber, S. Kaplan, M. Teicher, B. Tsaban and U. Vishne,
   {\it Length-based conjugacy search in the Braid group}. Preprint. Available at
   {\tt www.arxiv.org/math.GR/0209267.}



\bibitem{garside} F. A. Garside, {\it The braid group and other groups}. Quart. J.
Math. Oxford {\bf 20} (1969), 235-154.

\bibitem{rainer} D. Hofheinz and R. Steinwandt, {\it A Practical Attack on Some Braid Group Based
Cryptographic Primitives}. Accepted for presentation at the International Workshop on Practice and
Theory in Public Key Cryptography - PKC 2003.

\bibitem{H} J. Hughes, {\it A Linear Algebraic Attack on the AAFG1 Braid Group Cryptosystem}.
The 7th Australasian Conference on Information Security and Privacy ACISP 2002,
Lecture Notes in Computer Science, {\bf 2384}, 176--189, Springer-Verlag, New York  2002.

\bibitem{HT} J. Hughes and A. Tannenbaum, {\it Length-Based Attacks for Certain Group Based Encryption Rewriting Systems},
Workshop SECI02 Securit\'e de la Communication sur Intenet, September 2002, Tunis, Tunisia.


\bibitem{leelee} S. J. Lee and E. Lee, {\it Potential weaknesses
of the commutator key agreement protocol based on braid groups}. L.R. Knudsen (Ed.):
EUROCRYPT 2002, LNCS 2332, pp. 14-28, 2002.

\bibitem{MK} K. Murasugi and B. Kurpita, ``A Study of Braids", Kluwer, Dordrecht 1999.

\bibitem{T} W. P. Thurston, Braid Groups, Chapter 9 of ``Word processing in groups",
D. B. A. Epstein, J. W. Cannon, D. F. Holt, S. V. F. Levy, M. S. Paterson and
W. P. Thurston, Jones and Bartlett Publishers, Boston, MA, 1992.

\end{thebibliography}
\end{document}